# Metric transformations under collapsing of Riemannian manifolds

Bennett Chow
University of California, San Diego

David Glickenstein
University of California, San Diego

Peng Lu
University of Oregon

## 1 Introduction

Gromov-Hausdorff convergence is an important tool in comparison Riemannian geometry. Given a sequence of Riemannian manifolds of dimension $n$ with Ricci curvature bounded from below, Gromov's precompactness theorem says that a subsequence will converge in the pointed Gromov-Hausdorff topology to a length space [G-99, Section 5A]. If the sequence has bounded sectional curvature, then the limit will be an Alexandrov space [BGP-92]. The general idea in using Gromov-Hausdorff convergence is to construct sequences and to inspect their limits closely (see [P-97], for example). A sequence is called collapsing if its Gromov-Hausdorff limit has a lower Hausdorff dimension. It is known that interesting phenomena occur when a sequence of Riemannian manifolds collapses [CG-86], [CG-90], [Fu-90].

In this paper we combine the ideas of isometric group actions and collapsing of Riemannian manifolds to explore metric transformations. More precisely, for several Riemannian manifolds $(M, g)$ with isometric group actions $G$, we construct a sequence of Riemannian manifolds $(\widetilde{M}_i, \widetilde{g}_i)$ of higher dimension which collapse to the original manifold $M$ with a new metric $h$. $G$ also acts on $(M, h)$ by isometries. The map $g \mapsto h$ is the metric transformation in the title of this paper. Such changes in metrics were previously considered by Cheeger [C-72] to construct manifolds with nonnegative curvature. A new observation in this paper is that when $g$ is the hyperbolic metric on $\mathbb{H}^2$ and $G = S^1$, the transformed metric $h$ is Hamilton's cigar soliton metric studied in Ricci flow [H-88], [H-95] (see also [W-91, p. 315]). We also observe that when $g$ is the exploding soliton (see the appendix for its definition), the transformed metric $h$ is the 2-sphere. When $g$ is the standard metric on $S^3$ and $G = S^1$, the transformed metrics $h$ are the Berger metrics on $S^3$ [CE-75].

In section 2 we construct examples of sequences $(\widetilde{M}_i, \widetilde{g}_i)$ and show that the collapsing limit $(M, h)$ of $(\widetilde{M}_i, \widetilde{g}_i)$ is a quotient of a Riemannian manifold $(\widetilde{M}, \widetilde{g})$ by a Lie group $\widetilde{G}$. In section 3 we calculate the transformed metric $h$ for the examples in section 2. In section 4 we compute the metrics on quotients $(\widetilde{M}, \widetilde{g})/\widetilde{G}$ for other choices of $(\widetilde{M}, \widetilde{g})$ and $\widetilde{G}$ which give metric transformations.

Throughout this paper we shall use GH to denote Gromov-Hausdorff.

**Acknowledgement 1** *We would like to thank the National Center for Theoretical Sciences at Hsinchu, Taiwan for providing partial support during the summer of 2002 where the main portion of this work was completed. We also thank Cumrun Vafa for discussions during January of 2003 concerning Witten's black hole (Hamilton's cigar soliton). We are grateful to Kris Tapp for bringing to our attention the work of Cheeger [C-72].*

## 2 Examples of collapsing sequences

In this section we give some examples of collapsing sequences and describe their limits. There are two pictures of collapse we will consider: (A) collapse via quotients of finite groups with orders tending to infinity and (B) collapse as a fibration of a manifold over an Alexandrov space (which can sometimes be viewed as a limit of case (A)). In the first subsection we present examples of (A) and in the second subsection we present an example of (B).

### 2.1 Collapse via quotients of finite groups

Let $S^1(r)$ be the circle of radius $r$, i.e. $S^1(r) \doteq \mathbb{R}/(2\pi r \mathbb{Z})$. Let $(P^2, g)$ be a rotationally symmetric surface (compact or noncompact) with an isometric action of $S^1$ by rotations. Some examples of $P$ are the 2-sphere, hyperbolic disk, and Hamilton's cigar soliton metric. Given $\tau \in S^1 \doteq \mathbb{R}/(2\pi\mathbb{Z})$ and $x \in P$, we let $\tau(x) \in P$ denote the image of $x$ under the rotation by an angle $\tau$.

**Example 1**. Define a $\mathbb{Z}_p$-action on the product $P \times S^1(r)$ by the following: for $q \in \mathbb{Z}_p \doteq \mathbb{Z}/(p\mathbb{Z})$

$$x \mapsto \frac{2\pi q}{p}(x), \qquad y \mapsto y + \frac{2\pi r q}{p}.$$

The GH limit of the sequence $\{[P \times S^1(r)]/\mathbb{Z}_{p_i}\}_{i=1}^{+\infty}$, where $p_i \to +\infty$, is isometric to the quotient $[P \times S^1(r)]/S^1$ where the action of $\tau \in S^1$ is defined by

$$x \mapsto \tau(x), \qquad y \mapsto y + r\tau.$$

**Example 2**. Let $m_1$ and $m_2$ be two natural numbers and $\kappa \doteq m_1/m_2$. We define a $\mathbb{Z}_p$-action on $P \times S^1(r)$ with slope $\kappa$ by the following: for $q \in \mathbb{Z}_p$

$$x \mapsto \frac{2\pi m_1 q}{p}(x), \qquad y \mapsto y + \frac{2\pi r m_2 q}{p}.$$

The GH limit of the sequence $\{[P \times S^1(r)]/_\kappa \mathbb{Z}_{p_i}\}_{i=1}^{+\infty}$,[1] where $p_i \to +\infty$, is isometric to the quotient $[P \times S^1(r)]/_\kappa S^1$ with the action of $\tau \in S^1$ defined by

$$x \mapsto [m_1 \tau](x), \qquad y \mapsto y + r m_2 \tau.$$

---

[1] We use the subscript $\kappa$ to emphasize the dependence of the action on the slope $\kappa$.

A simpler 2-dimensional version of Example 2 is:

**Example 3**. Define a $\mathbb{Z}_p$-action on $S^1(r_1) \times S^1(r_2)$ with slope $\kappa = m_1/m_2$ by the following: for $q \in \mathbb{Z}_p$

$$x \mapsto x + \frac{2\pi r_1 m_1 q}{p}, \qquad y \mapsto y + \frac{2\pi r_2 m_2 q}{p}$$

The GH limit of the sequence $\{[S^1(r_1) \times S^1(r_2)]/_\kappa \mathbb{Z}_{p_i}\}_{i=1}^{+\infty}$, where $p_i \to \infty$, is isometric to the quotient $[S^1(r_1) \times S^1(r_2)]/_\kappa S^1$ with the action of $\tau \in S^1$ defined by

$$x \mapsto x + m_1 r_1 \tau, \qquad y \mapsto y + m_2 r_2 \tau.$$

## 2.2 Collapse as a fibration of manifolds

In this section we look at the Berger spheres. First consider the Berger spheres from the point of view given in [CE-75, 3.35]. The Lie algebra of the Lie group $SU(2) = S^3$ has a basis of left invariant vector fields $F_1, F_2, F_3$ which satisfy

$$[F_1, F_2] = 2F_3, \quad [F_2, F_3] = 2F_1, \quad [F_3, F_1] = 2F_2.$$

The three vector fields are orthonormal in the standard metric on $S^3$. Let $F_4$ be the left invariant vector field on $S^1$ and let $\widetilde{g}$ be the standard product metric on $S^3 \times S^1$. For each $\xi \in S^1$ the vector field $\cos\xi \cdot F_1 + \sin\xi \cdot F_4$ defines an $S^1$-isometric action on $S^3 \times S^1$. We denote the quotient by $(S^3 \times S^1)/_\xi S^1$. If $\xi \neq 0$ or $\pi$ the quotient is diffeomorphic to $S^3$ and we get a Berger sphere metric on $S^3$. In particular, if $\xi = \pi/2$ or $3\pi/2$ we get the standard metric on $S^3$. We can think of $\xi$ as slope of the action. As $\xi \to 0$ or $\pi$ the family of Berger spheres $(S^3 \times S^1)/_\xi S^1$ converges (collapses) to a standard sphere $S^2$ in the GH topology. As $\xi \to \pi/2$ or $3\pi/2$ the family of Berger spheres $(S^3 \times S^1)/_\xi S^1$ converges to a standard $S^3$ in the GH topology.

We can also look at the Berger spheres via the Hopf fibration. Define left invariant metrics on the Lie group $SU(2)$ such that $F_1, F_2, F_3$ are orthogonal but not necessarily orthonormal. These are the left invariant metrics as described in Milnor [M-76] and, in the context of Ricci flow, in Isenberg-Jackson [IJ-92]. Following their treatment we consider the metrics

$$A\omega^1 \otimes \omega^1 + B\omega^2 \otimes \omega^2 + C\omega^3 \otimes \omega^3$$

where $\omega^i$ is dual to $F_i$ for $i = 1, 2, 3$. Collapsing to $S^2(2B)$, the sphere of radius $2B$, corresponds to $A \to 0^+$ and $B = C$. We can see this by considering the Hopf fibration $H: S^3 \to S^2$ defined by $H(z,w) = [z,w]$ where $(z,w) \in \mathbb{C}^2$ such that $|z|^2 + |w|^2 = 1$ and $[z,w]$ are homogeneous coordinates on $S^2 = \mathbb{CP}^1$. The vector field $F_1$ is tangent to the fiber, which is homeomorphic to $S^1$. A simple calculation shows that if $B = C$ then the Hopf fibration is a Riemannian submersion from the Berger sphere metrics to the standard metric on $S^2(2B)$. When the length of the fiber is less than $\varepsilon$, $H$ is an $\varepsilon$-GH approximation and it follows that as $A$ goes to 0 the Berger spheres converge to $S^2(2B)$.

# 3 Examples of metric transformations

For the rest of this paper we will use the following conventions.

(i) The coordinate on the circle $S^1(r)$ is $0 \leq s < 2\pi$ and the metric on $S^1(r)$ is $r^2 ds^2$.

(ii) The coordinates on the rotationally symmetric surface $(P^2, g)$ are $(\rho, \theta)$ with $\rho \in I \subset \mathbb{R}$, an interval, and $\theta \in [0, 2\pi)$. The metric is

$$g = d\rho^2 + f(\rho)^2 d\theta^2.$$

(iii) The action of $\tau \in S^1 \doteq \mathbb{R}/(2\pi\mathbb{Z})$ on $P$ is defined by

$$\tau(\rho, \theta) = (\rho, \theta + \tau). \tag{3.1}$$

## 3.1 Metrics on $S^1$-quotients from subsection 2.1

**Example 1**. (The quotient $[(P^2, g) \times S^1(r)]/S^1$) The $S^1$-action on $(P^2, g) \times S^1(r)$ is defined by

$$\rho \mapsto \rho, \qquad \theta \mapsto \theta + \tau, \qquad s \mapsto s + \tau$$

for $\tau \in S^1$. The map $q : P^2 \to [P^2 \times S^1(r)]/S^1$, where

$$q(\rho, \theta) = [\rho, \theta, 0],$$

is a diffeomorphism and gives coordinates $(\rho, \theta)$ on the quotient space $[(P^2, g) \times S^1(r)]/S^1$.

**Proposition 2** *The quotient metric $h$ on $[(P^2, g) \times S^1(r)]/S^1$ is*

$$h = d\rho^2 + \frac{r^2 f(\rho)^2}{f(\rho)^2 + r^2} d\theta^2.$$

**Proof.** This computation is a special case of [C-72, p.624]. For the convenience of the reader we give the calculation below. The metric on $P^2 \times S^1(r)$ is

$$g_{\text{prod}} = d\rho^2 + f(\rho)^2 d\theta^2 + r^2 ds^2.$$

The unit vector field tangent to the orbits of the $S^1$ action is

$$W = \frac{1}{\sqrt{f(\rho)^2 + r^2}} \left( \frac{\partial}{\partial \theta} + \frac{\partial}{\partial s} \right).$$

We have a submersion $p : P^2 \times S^1(r) \to [P^2 \times S^1(r)]/S^1$, where $p(\rho, \theta, s) = [\rho, \theta, s]$. The differential $p_* : T[P^2 \times S^1(r)] \to T[P^2 \times S^1(r)]/S^1$ is given by

$$p_*(V) = V - g_{\text{prod}}(V, W)W.$$

We have

$$p_*\left(\frac{\partial}{\partial \rho}\right) = \frac{\partial}{\partial \rho},$$

$$p_*\left(\frac{\partial}{\partial \theta}\right) = \frac{r^2}{f(\rho)^2 + r^2}\frac{\partial}{\partial \theta} - \frac{f(\rho)^2}{f(\rho)^2 + r^2}\frac{\partial}{\partial s},$$

$$p_*\left(\frac{\partial}{\partial s}\right) = -\frac{r^2}{f(\rho)^2 + r^2}\frac{\partial}{\partial \theta} + \frac{f(\rho)^2}{f(\rho)^2 + r^2}\frac{\partial}{\partial s}.$$

Using $h(X,Y) = g_{\text{prod}}(p_*(X), p_*(Y))$, we find

$$h\left(\frac{\partial}{\partial \rho}, \frac{\partial}{\partial \rho}\right) = 1, \quad h\left(\frac{\partial}{\partial \rho}, \frac{\partial}{\partial \theta}\right) = 0, \quad h\left(\frac{\partial}{\partial \theta}, \frac{\partial}{\partial \theta}\right) = \frac{r^2 f(\rho)^2}{f(\rho)^2 + r^2}.$$

The proposition is proved. ∎

The same method applies to the following two examples.

**Example 2**. (The quotient $[(P^2, g) \times S^1(r)]/_\kappa S^1$) Let $m_1, m_2$ be two natural numbers and $\kappa \doteq m_1/m_2$. The $S^1$-action with slope $\kappa$ on $(P^2, g) \times S^1(r)$ is defined by

$$\rho \mapsto \rho, \quad \theta \to \theta + m_1\tau, \quad s \mapsto s + m_2\tau,$$

for $\tau \in S^1$. The quotient $[(P^2, g) \times S^1(r)]/_\kappa S^1$ is a surface diffeomorphic to $P^2$ and the quotient metric $h_\kappa$ is

$$h_\kappa = d\rho^2 + \frac{r^2 f(\rho)^2}{\kappa^2 f(\rho)^2 + r^2} d\theta^2,$$

where the coordinates $(\rho, \theta)$ are as in Example 1.

**Example 3**. (The quotient $[S^1(r_1) \times S^1(r_2)]/_\kappa S^1$) The product metric on $S^1(r_1) \times S^1(r_2)$ is $g_{\text{prod}} = r_1^2 ds_1^2 + r_2^2 ds_2^2$ and the $S^1$-action with slope $\kappa = m_1/m_2$ is defined by

$$s_1 \mapsto s_1 + m_1\tau, \quad s_2 \mapsto s_2 + m_2\tau,$$

for $\tau \in S^1$. The map $q : S^1(r_1) \to [S^1(r_1) \times S^1(r_2)]/_\kappa S^1$, defined by $q(s_1) = [s_1, 0]$, gives a coordinate $s_1 \in [0, 2\pi)$ on $[S^1(r_1) \times S^1(r_2)]/_\kappa S^1$. The quotient $[S^1(r_1) \times S^1(r_2)]/_\kappa S^1$ is isometric to the circle $S^1(\sqrt{\frac{r_1^2 r_2^2}{\kappa^2 r_1^2 + r_2^2}})$, i.e., the quotient metric is

$$h = \frac{r_1^2 r_2^2}{\kappa^2 r_1^2 + r_2^2} ds_1^2.$$

## 3.2 Metric transformation on surfaces admitting an $S^1$ action

Motivated by the formula in example 2 in section 3.1 we give

**Definition 3** *Let $r > 0$ and $\mathcal{G}_{S^1}^P$ be the set of rotationally symmetric metrics on the surface $P$. We define a metric transformation $\Upsilon_{r,\kappa} : \mathcal{G}_{S^1}^P \to \mathcal{G}_{S^1}^P$ by*

$$\Upsilon_{r,\kappa}(d\rho^2 + f(\rho)^2 d\theta^2) = d\rho^2 + \frac{r^2 f(\rho)^2}{\kappa^2 f(\rho)^2 + r^2} d\theta^2.$$

*We denote $\Upsilon_{1,1}$ by $\Upsilon$.*

We calculate a few special examples of $\Upsilon(g)$.

**1**. If $P = \mathbb{R}^1 \times S^1$, $g = d\rho^2 + d\theta^2$, then $\Upsilon(g) = d\rho^2 + \frac{1}{2}d\theta^2$. That is, the flat cylinder of radius 1 transforms to the flat cylinder of radius $1/\sqrt{2}$.

**2**. If $P = \mathbb{R}^2$, $g_h \doteq d\rho^2 + \sinh^2 \rho \, d\theta^2$, then

$$\Upsilon(g_h) = d\rho^2 + \tanh^2 \rho \, d\theta^2.$$

Note that $g_h$ is the hyperbolic metric and $\Upsilon(g_h)$ is Hamilton's cigar soliton metric [H-88], [H-95] (also known as the Witten black hole [W-91]). This shows how the cigar soliton metric is the GH limit of the sequence $\{[P \times S^1(1)]/\mathbb{Z}_{p_i}\}_{i=1}^{+\infty}$ as $p_i \to +\infty$.

**3**. If $P = S^2$, $g = d\rho^2 + \tan^2 \rho \, d\theta^2$, then

$$\Upsilon(g) = d\rho^2 + \sin^2 \rho \, d\theta^2,$$

which is the standard metric on the upper hemisphere of $S^2(1)$. The metric $g$ satisfies the equation $\Delta \ln(-R) + R = 0$ and its scalar curvature is $R = -4\sec^2 \rho$, which tends to $-\infty$ as $\rho \to \pi/2$. The incomplete metric $g$ is the *exploding soliton* described in the appendix. This example seems in some sense dual to the previous example.

We find it remarkable that the hyperbolic disk $\mathbb{H}^2$ is transformed to Hamilton's cigar soliton $\Sigma^2$ and that the exploding soliton is transformed to the 2-sphere. In particular, for simply connected two dimensional rotationally symmetric metrics, the complete *trivial expanding* Ricci soliton is transformed into the complete noncompact steady Ricci soliton with *positive* curvature. The maximal noncompact steady Ricci soliton with *negative curvature* is transformed to the *trivial shrinking* compact Ricci soliton. Is there a more general reason for this? It would also be interesting to find other special metrics which are transformed to special metrics.

## 4 Transformations of metrics with symmetries

A transformation of metrics more general than that in subsection 3.2 was given by Cheeger in [C-72], where he considered $(M \times G_1)/G$, with $G$ a closed subgroup of a Lie group $G_1$. In subsection 4.1 we give some special cases of this transformation. In subsection 4.2 we give a generalization.

## 4.1 Transformation via a diagonal $S^1$-action with slope $\kappa$

Let $M^n$ be a manifold with an isometric $S^1$ group action $\mathcal{I}: S^1 \times M \to M$. The $S^1$ action induces a Killing vector field $K$ on $(M, g)$ by

$$K(x) \doteqdot \left.\frac{\partial}{\partial \tau}\right|_{\tau=0} \mathcal{I}(\tau, x).$$

Let $r > 0$ and $m_1, m_2$ be two natural numbers. We consider an $S^1$-action $\mathcal{J}$ on $\left(M \times S^1(r), g + r^2 ds^2\right)$ with slope $\kappa = m_1/m_2$ defined by

$$\mathcal{J}(\tau, x, s) \doteqdot (\mathcal{I}(m_1 \tau, x), s + m_2 \tau).$$

We will denote the quotient space by $\left(M \times S^1(r)\right)/_\kappa \mathcal{J}$.

By [C-72, p. 624], we have

**Proposition 4** *Let $M$ be a manifold admitting an $S^1$ action and $\mathcal{G}_{S^1}^M$ be the set of Riemannian metrics $g$ such that $S^1$ acts on $(M, g)$ by isometries. Then the induced metric $h_{r,\kappa}$ on the quotient $\left(M \times S^1(r)\right)/_\kappa \mathcal{J}$ defines a metric transformation*

$$\Upsilon_{r,\kappa}: \mathcal{G}_{S^1}^M \to \mathcal{G}_{S^1}^M,$$

$$\Upsilon_{r,\kappa}(g) = g - \frac{\kappa^2}{\kappa^2 |K|^2 + r^2} K^* \otimes K^*,$$

*where $K^*$ is the 1-form dual to $K$ (using the metric $g$).*

Note the above metric transformation is well-defined for Riemannian manifolds admitting a Killing vector field.

As an application of proposition 4, we show that the standard metric $g_{std}$ on $S^3(1)$ can be transformed to the Berger metrics. We adopt the notation of section 2.2. Consider the isometric $S^1$ action on $(S^3(1), g_{std})$ whose Killing field is $K = F_1$. Choose the metric $g$ in proposition 4 to be

$$g_{std} = \omega^1 \otimes \omega^1 + \omega^2 \otimes \omega^2 + \omega^3 \otimes \omega^3,$$

$|K| = 1$ and $K^* = \omega^1$. Hence when choosing $r = 1$ the transformed metric is

$$\Upsilon_{1,\kappa}(g_0) = \frac{1}{\kappa^2 + 1} \omega^1 \otimes \omega^1 + \omega^2 \otimes \omega^2 + \omega^3 \otimes \omega^3,$$

which is a Berger metric.

## 4.2 Transformations from products

Let $(M_\alpha^{n_\alpha}, g_\alpha), \alpha = 1, 2$ be two Riemannian manifolds, each admitting $m$ Killing vector fields (which we think of as coming from the action of the same local group on both manifolds). Denote the Killing vector fields on $M_\alpha$ by $K_{\alpha i}$ for $i = 1, \ldots, m$. Endow $M_1 \times M_2$ with the product metric $g_{\text{prod}} = g_1 + g_2$.

Define the subspace $H$ of the tangent space $T(M_1 \times M_2)$ by

$$H = \text{span}\{(K_{1i}, K_{2i}) : i = 1, \ldots, m\}.$$

Suppose there is an immersed submanifold $q : N \to M_1 \times M_2$ of dimension $n_1 + n_2 - m$ such that for any $(x, y) \in N$, $H$ is transverse to $N$ in $M_1 \times M_2$. We shall define a metric $h_N$ on $N$.

Let $P_H : T(M_1 \times M_2) \to H$ be the orthonormal projection with respect to the metric $g_{\text{prod}}$. Suppose $X \in TN$. The projection of $q_* X \in T(M_1 \times M_2)$ onto the subspace $H^\intercal$, which is the space perpendicular to $H$ in $T(M_1 \times M_2)$, is given by

$$X^\top \doteq q_* X - P_H(q_* X).$$

Then for any $X, Y \in TN$, we define the metric $h_N$ on $N$ by

$$h_N(X, Y) = g_{\text{prod}}\left(X^\top, Y^\top\right). \tag{4.4}$$

We call the map $((M_1, g_1), (M_2, g_2)) \to (N, h_N)$ the metric transformation. Note that this metric transformation is not the same as in the previous sections in the sense that two Riemannian manifolds with isometric group actions are transformed to a new Riemannian manifold.

As a special case we consider the following situation. Let $\mathbb{G}^m$ be a Lie group acting isometrically on $(M_\alpha, g_\alpha)$

$$\mathcal{I}_\alpha : \mathbb{G} \times M_\alpha \to M_\alpha.$$

Define an action

$$\mathcal{J} : \mathbb{G} \times (M_1 \times M_2) \to M_1 \times M_2$$

by

$$\mathcal{J}(\tau, x, y) \doteq (\mathcal{I}_1(\tau, x), \mathcal{I}_2(\tau, y)).$$

$\mathcal{J}$ acts isometrically on $M_1 \times M_2$ with the product metric $g_{\text{prod}}$.

Let $p : M_1 \times M_2 \to (M_1 \times M_2)/\mathcal{J}$ be the projection map and suppose $q : N \to M_1 \times M_2$ is an immersed submanifold transverse to the orbits of $\mathbb{G}$. Then the map $p \circ q : N \to (M_1 \times M_2)/\mathcal{J}$ gives coordinates on $(M_1 \times M_2)/\mathcal{J}$. It is clear that the metric $h_N$ defined above is the pullback by $p \circ q$ of the quotient metric on $(M_1 \times M_2)/\mathcal{J}$.

## 5 Appendix: The cigar and exploding solitons

Hamilton's **cigar soliton** $\Sigma^2$ (see below for why it is a Ricci soliton) is the complete Riemannian surface with underlying manifold $\mathbb{R}^2$ and metric

$$g_\Sigma = \frac{dx^2 + dy^2}{1 + x^2 + y^2}.$$

In polar coordinates, this is $g_\Sigma = \frac{dr + r^2 d\theta^2}{1 + r^2}$, which exhibits the fact that the cigar is asymptotic to a cylinder of radius 1. Making the change of variables $\rho = \sinh^{-1} r = \ln\left(\sqrt{1+r^2} + r\right)$, we have $g_\Sigma = d\rho^2 + \tanh^2 \rho\, d\theta^2$. The Gauss curvature is given by $K = 2\cosh^{-2}\rho$, which decays exponentially in $\rho$.

Consider a rotationally symmetric surface with a metric of the form

$$g = d\rho^2 + f(\rho)^2 d\theta^2,$$

where $f(\rho)$ is some positive function. The Gauss curvature is given by $K = -f''(\rho)/f(\rho)$. Recall that $g$ is a steady gradient Ricci soliton if there exists a function $\phi$ such that

$$Kg_{ij} = \nabla_i \nabla_j \phi. \tag{A1}$$

Assume that $\phi$ is a radial function $\phi = \phi(\rho)$. In terms of the orthonormal frame, (A1) implies

$$-\frac{f''}{f} = \phi'' = \frac{f'}{f}\phi'.$$

We deduce

$$f' + Af^2 = B$$

for some constants $A$ and $B$. Now assume that $g$ extends to a smooth metric at $\rho = 0$, that is, $f(0) = 0$ and $f'(0) = 1$. This implies $B = 1$. Solving the ODE gives

**1**. If $A = a^2 > 0$, then $f(\rho) = \frac{1}{a}\tanh(a\rho)$. Thus

$$g = d\rho^2 + \frac{1}{a^2}\tanh^2(a\rho)\, d\theta^2, \quad \text{where } a > 0,$$

which is a constant multiple of the cigar soliton. Note that the potential function is given by $\phi(\rho) = 2\ln(\cosh(a\rho))$.

**2**. If $A = -a^2 < 0$, then $f(\rho) = \frac{1}{a}\tan(a\rho)$. In this case

$$g = d\rho^2 + \frac{1}{a^2}\tan^2(a\rho)\, d\theta^2, \quad \text{where } a > 0.$$

We call this the **exploding soliton** since $f(\rho)$ explodes to infinity as $a\rho \to \pi/2$. Note that the potential function is given by $\phi(\rho) = 2\ln(\cos(a\rho))$. This metric appears in [DVV-92, p. 295].